\documentclass[a4paper,11pt]{amsart}
\usepackage{amssymb}

\theoremstyle{plain}
\newtheorem{theorem}[subsection]{Theorem}

\newtheorem{lemma}[subsection]{Lemma}
\newtheorem{proposition}[subsection]{Proposition}

\theoremstyle{definition}

\newtheorem{example}[subsection]{Example}
\theoremstyle{remark}
\newtheorem{remark}[subsection]{Remark}
\numberwithin{equation}{section}

\newcommand{\kc}{{\mathcal C}}

\newcommand{\ko}{{\mathcal O}}

\newcommand{\kv}{{\mathcal V}}

\newcommand{\IC}{{\mathbb C}}
\newcommand{\IH}{{\mathbb H}}

\newcommand{\IN}{{\mathbb N}}
\newcommand{\IP}{{\mathbb P}}
\newcommand{\IQ}{{\mathbb Q}}

\newcommand{\IZ}{{\mathbb Z}}

\newcommand{\gothf}{{\mathfrak f}}
\newcommand{\gothg}{{\mathfrak g}}
\newcommand{\gothl}{{\mathfrak l}}
\newcommand{\gothp}{{\mathfrak p}}

\newcommand{\gothr}{{\mathfrak r}}

\newcommand{\id}{{\rm id}}
\newcommand{\tensor}{\otimes}

\newcommand{\Hilb}{{\rm Hilb}}
\newcommand{\lra}{\longrightarrow}
\newcommand{\xra}{\xrightarrow}
\newcommand{\isom}{\cong}

\newcommand{\epimorph}{\twoheadrightarrow}
\DeclareMathOperator{\Frobint}{T}

\DeclareMathOperator{\vol}{vol}

\DeclareMathOperator{\End}{End}
\DeclareMathOperator{\ad}{ad}
\DeclareMathOperator{\glf}{\gothg\gothl\gothf}
\newcommand{\eps}{\varepsilon}
\DeclareMathOperator{\sgn}{sgn}
\DeclareMathOperator{\vacuum}{\mathbf{1}}

\newcommand{\bijection}{\cong}
\newcommand{\half}{\frac{1}{2}}
\newcommand{\<}{\langle}
\renewcommand{\>}{\rangle}
\newcommand{\orb}[2]{\<#1\>\backslash #2}
\newcommand{\dprime}{{\prime\prime}}

\newcommand{\Nakop}{\gothp}
\begin{document}
\title[Cup Product]{The cup product of the Hilbert scheme for K3 surfaces}
\author{Manfred Lehn \and Christoph Sorger}
\date{December 15, 2000}

\begin{abstract}To any graded Frobenius algebra $A$ we associate a sequence of
graded Frobenius algebras $A^{[n]}$ so that there is canonical isomorphism
of rings $(H^*(X;\IQ)[2])^{[n]}\isom H^*(X^{[n]};\IQ)[2n]$ for the
Hilbert scheme $X^{[n]}$ of generalised $n$-tuples of any smooth projective
surface $X$ with trivial canonical bundle.
\end{abstract}

\address{
Manfred Lehn\\
Mathematisches Institut
der Universit\"at zu K\"oln\\
Weyertal 86-90\\
D-50931 K\"oln, Germany}
\email{manfred.lehn@math.uni-koeln.de}
\address{
Christoph Sorger\\
Math\'{e}matiques (UMR 6629 du CNRS)\\
Universit\'{e} de Nantes\\
2, Rue de la Houssini\`{e}re\\
BP 92208\\
F-44322 Nantes Cedex 03, France}
\email{christoph.sorger@math.univ-nantes.fr}

\maketitle


\section{Introduction}

Let $X$ be a smooth projective surface and let $X^{[n]}$ denote the
Hilbert scheme of generalized $n$-tuples for each non-negative integer 
$n$. The topological or geometric data of $X^{[n]}$ can be expressed in 
terms of the corresponding data of the surface itself intricately interwoven
with the combinatorics of the symmetric group. Instances of this dogma
are G\"ottsche's Formula for the Betti numbers and Nakajima's and Grojnowski's
lift of this formula to the level of vector spaces with non-degenerate
pairings:
$$\bigoplus_{n\geq 0} H^*(X^{[n]};\IQ)\xra{\ \isom\ }S^*(t^{-1}\IQ[t^{-1}]
\tensor H^*(X;\IQ))=\kv({H^*(X;\IQ)}),$$
where $\kv({H^*(X;\IQ)})$ stands for the bosonic vertex algebra modelled
on the vector space $H^*(X;\IQ)$ equipped with the cup product pairing.

This isomorphism allows one to reconstruct the cup-product pairing on 
$H^*(X^{[n]};\IQ)$, but not the cup product itself. The cup product
structure is the topic of this paper. A complete description of the
product in the framework of the vertex algebra calculus was given in
\cite{Lehn} for surfaces whose cohomology is generated
by classes that can be deformed to algebraic ones, and by Li, Qin and Wang
\cite{LiQinWang1} for arbitrary projective surfaces. However, these
descriptions are rather implicit in terms of multiplication operators or
explicit generators with implicitly given relations. 

Our intention here is to construct a sequence of endofunctors $A\mapsto
A^{[n]}$ on the category of Frobenius algebras, which when applied to the
cohomology ring of a projective surface yield the cohomology rings of its
Hilbert schemes. This generalises the authors' \cite{LehnSorger} and 
Vasserot's \cite{Vasserot} theorem that there is a canonical ring isomorphism
$\mbox{\rm gr}^F\kc(S_n)\lra H^*((\IC^2)^{[n]};\IZ)$, where $\kc(S_n)$ is the
centre of the group ring $\IZ[S_n]$, and $F$ is a natural filtration.
Though the affine plane is not even projective and hence is not quite 
covered by the formalism we are going to describe, it has the advantage 
that its topology is trivial and we may undisturbedly observe the combinatoric
effects. The cue to our proof of the theorem was provided by an observation
of Frenkel and Wang about the similarity of certain differential operators
in the work of Goulden \cite{Goulden} and the first named author \cite{Lehn}.

The paper is organized as follows: in the first part we associate to any graded
Frobenius algebra $A$ a sequence of algebras $A\{S_n\}$, which in some sense lie
in between the group algebras $A[S_n]$ and the wreath product algebras
$A^{\tensor n}S_n$.  Their invariant subrings $A^{[n]}:=(A\{S_n\})^{S_n}$ are
the desired new Frobenius algebras. We show that $\bigoplus_{n\geq 0}A^{[n]}
\isom \kv({A})$. In the second part we set up the necessary background about
the geometry of Hilbert schemes. In the third we show that if
$\kv(H^*(X;\IQ))$ is endowed with this ring structure then Nakajima's
isomorphism is an isomorphism of rings, if $X$ has trivial canonical divisor.

The presence of a non-trivial canonical divisor requires a deformation 
of this product due to a curious correction term in \cite[Theorem 3.10.2]{Lehn}.

Lothar G\"ottsche informed us that Barbara Fantechi and he have shown that our
construction $A\to A^{[n]}$ identifies to a slightly modified version of
Chen's and Ruan's orbifold cohomology
construction in the case of the $n$-fold symmetric product of $X$ so that our
main theorem implies Y.\ Ruan's conjecture 6.3 in \cite{Ruan}
that the orbifold cohomology ring is isomorphic to the cohomology ring of
the Hilbert scheme of $n$-points.


\section{The algebraic model}


\subsection{} In this section we will construct a sequence of endofunctors
$A\to A^{[n]}$ in the category of graded Frobenius algebras. 

For our purposes a graded Frobenius algebra of degree $d$ is a finite dimensional
graded vector space $A=\bigoplus_{i=-d}^{d} A^i$ with a graded commutative
and associative multiplication $A\tensor A\to A$ of degree $d$ and unit
element 1 (necessarily of degree $-d$) together with a linear form $\Frobint:A\to\IQ$
of degree $-d$ such that the induced symmetric bilinear form
$\langle a,b\rangle\,:=\,\Frobint(ab)$ is non-degenerate (and of degree 0). It 
follows from $1\cdot 1=1$ that $d$ must be an even number. The degree 
conventions are chosen in such a way that $A$ is centred around degree 0.
The degree of an element $a$ will be denoted by $|a|$.

In the applications, $A$ will be the shifted cohomology ring $H^*(X;\IQ)[d]$ of
a compact complex manifold $X$ of even dimension $d$.

The tensor product $A^{\tensor n}$ is again a graded Frobenius algebra 
of degree $nd$ with product
$$(a_1\tensor\ldots\tensor a_n)\cdot(b_{1}\tensor\ldots\tensor b_{n}):=
\varepsilon(\underline{a},\underline{b})(a_1b_1)\tensor\ldots\tensor 
(a_nb_n),$$
where $\varepsilon(\underline{a},\underline{b})$ is the sign resulting from
reordering the $a$'s and $b$'s. The integral is given by
$\Frobint(a_1\tensor\ldots\tensor a_n):=\Frobint(a_1)\cdot\ldots\cdot
\Frobint(a_n$.)

The symmetric group $S_n$ acts on the $n$-fold tensor product $A^{\tensor n}$ 
as
$$\pi(a_1\tensor\ldots\tensor a_n)=
\varepsilon(\pi,\underline{a})a_{\pi^{-1}(1)}\tensor\ldots\tensor a_{\pi^{-1}(n)},$$
where $\varepsilon(\pi,\underline{a})=(-1)^{\sum_{i<j, 
\pi(j)<\pi(i)}|a_i|\cdot|a_j|}$ is the sign introduced by interchanging
the $a_i$'s.

It will be useful to extend the definition of the $n$-fold 
tensor product to allow arbitrary unordered finite indexing sets.
Therefore let $I$ be a finite set with $n$ elements and let $\{A_i\}_{i\in I}$ 
be a family of copies of a graded Frobenius algebra $A$ indexed by $I$. Let 
$[n]$ denote the set $\{1,\ldots, n\}$. Then we define
$$A^{\tensor I}:=\left(\bigoplus_{f:[n]\xra{\bijection} I} A_{f(1)}\tensor 
\ldots\tensor A_{f(n)}\right)\Big/{S_n}.$$
Here the direct sum inherits a Frobenius algebra structure by componentwise
defined operations, the symmetric group acts on this direct sum via the
induced operation on the set of bijections $[n]\xra{\bijection} I$, and we
take the quotient by the action of $S_n$. Thus for each bijection $f:[n]\to I$
there is a canonical isomorphism $A^{\tensor n}\xra{\,\isom\,}A^{\tensor I}$.

Let $n_1,\ldots,n_k$ be natural numbers and let $n:=n_1+\ldots+n_k$. Consider
the ring homomorphism $\varphi_{n_\bullet,k}:A^{\tensor n}\to A^{\tensor k}$ 
which sends
$$a_1\tensor \ldots \tensor a_n\mapsto(a_1\cdots a_{n_1})\tensor \ldots 
\tensor(a_{n_1+\ldots+n_{k-1}+1}\cdots a_{n}).
$$
As before we want to generalise this to get a graded ring homomorphism 
$\varphi^*:A^{\tensor I}\to A^{\tensor J}$ for any surjective map
$\varphi:I\to J$ of finite sets of cardinality $n:=|I|$ and
$k:=|J|$. Choose a bijection $g:[k]\to J$ and let
$n_i:=|\varphi^{-1}(g(i))|$. Then there is a bijection $f:[n]\to I$ such 
that for each $i\in [k]$ one has 
$$\varphi^{-1}(g(i))=\{f(n_1+\ldots+n_{i-1}+1),\ldots,f(n_1+\ldots+n_i)\}.$$
The composition 
$$\varphi^*:A^{\tensor I}\xra{f^{-1}}A^{\tensor n}\xra{\varphi_{n_\bullet,k}} 
A^{\tensor k}\xra{\ g\ } A^{\tensor J}$$
is well-defined independently of the choices of $f$ and $g$.
If $\varphi:I\to J$ and $\psi:J\to K$ are surjections then 
$\psi^*\circ \varphi^*=(\psi\circ\varphi)^*$.

Finally, for any surjection $\varphi:I\to J$ let 
$$\varphi_*:A^{\tensor J}\to A^{\tensor I}$$
be the linear map adjoint to $\varphi^*$ with respect to the bilinear forms on
$A^{\tensor I}$ and $A^{\tensor J}$. Then $\varphi^*$ and $\varphi_*$ are both
homogeneous maps of degree $d(|I|-|J|)$. Whereas $\varphi^*$ is a ring
homomorphism, $\varphi_*$ is a module homomorphism with respect to $\varphi^*$, 
i.e.\ the `projection formula' holds:
$$\varphi_*(a\cdot \varphi^*(b))=\varphi_*(a)\cdot b$$
for all $b\in A^{\tensor I}$ and $a\in A^{\tensor J}$. 


\subsection{}
Consider the composite map
\begin{equation}\label{eq:simplestsurjection}
A\xra{\Delta_*}A\tensor A\lra A,
\end{equation}
where the second map is multiplication and $\Delta_*$ is the adjoint 
comultiplication. The image of $1$ under the composite linear map is called 
the Euler class of $A$ and is denoted by $e:=e(A)$. Note that if 
$A$ is connected, i.e. $\dim(A^{-d})=1$, then by duality we have 
$\dim(A^{d})=1$ as well. 
Hence there is a unique element $\vol\in A^{d}$ such that
$\Frobint(\vol)=1$. In this case $e(A)=\chi(A)\vol$, where 
$\chi(A):=\sum_i(-1)^i\dim(A^i)$ is the Euler-Poincar\'e-characteristic of $A$.  


\subsection{} For any numerical function $\nu:I\to \IN_0$ let
$$e^\nu:=\tensor_{i\in I}e^{\nu(i)}\ \in A^{\tensor I}.$$
For instance, if $\varphi:I\to J$ is a surjective map and
$\nu(j):=|\varphi^{-1}(j)|-1$ for
all $j\in J$, then
\begin{equation}
\varphi^*\varphi_*(a)=e^{\nu}\cdot a
\end{equation}
for all $a\in A^{\tensor J}$. More generally, let
\begin{equation}
\square:=\quad
\begin{array}{ccc}
I&\xra{\ \alpha }&K\\
\scriptstyle{\beta}\Big\downarrow\phantom{\scriptstyle{\beta}}
&&
\phantom{\scriptstyle{\gamma}}\Big\downarrow\scriptstyle{\gamma}\\
J&\xra{\ \delta\ }& L
\end{array}
\end{equation}
be a cocartesian diagram of finite sets and surjective maps. The associated
diagram
\begin{equation}
\begin{array}{ccc}
A^{\tensor I}&\xra{\ \alpha^*\ }&A^{\tensor K}\\
\scriptstyle{\beta_*}\Big\uparrow\phantom{\scriptstyle{\beta_*}}
&&
\phantom{\scriptstyle{\gamma_*}}\Big\uparrow\scriptstyle{\gamma_*}\\
A^{\tensor J}&\xra{\ \delta^*\ }& A^{\tensor L}
\end{array}
\end{equation}
of linear maps is not commutative. The deviation from commutativity is measured by the
following function: for any $\ell\in L$ let 
$$\nu_\square(\ell):=1-|\delta^{-1}(\ell)|-|\gamma^{-1}(\ell)|+
|(\delta\beta)^{-1}(\ell)|.$$
Then:

\begin{lemma}\label{lemma:diamond}--- $\alpha^*\beta_*(a)=
\gamma_*\left(e^{\nu_\square}\delta^*(a)\right)$ for all $a\in A^{\tensor L}$.
\end{lemma}

\begin{proof} It suffices to consider the case $|L|=1$. Let $G$ denote the
graph whose vertices consist of the sets  $J$ and $K$, and 
whose edges consist of all pairs $(\beta(i),i)$ and $(i,\alpha(i))$
for $i\in I$. The assumption that the diagram be cocartesian is equivalent to
the connectedness of the graph $G$. Every loop of this graph
gives rise to a contraction of type (\ref{eq:simplestsurjection}) and hence
introduces a factor $e$. The number of such loops is $1-(|J|+|K|)+|I|$.
\end{proof}

We leave the verification of the following observation to the reader (recall
that multiplication is a homogeneous map of degree $d$).

\begin{lemma}\label{edegree}--- Let $\nu:I\to \IN_0$ be a function on the finite set
$I$. Then $e^{\nu}\in A^{\tensor I}$ has degree
$|e^{\nu}|=2d\sum_{B\in I}\nu(B)-d|I|$.
\end{lemma}

\subsection{}\label{subsec:graphdefect} For any permutation $\pi\in S_n$ we call 
$$|\pi|:=\min\{m|\,\exists\text{ transpositions } \tau_1,\ldots,\tau_m 
\text{ such that }\pi=\tau_1\cdot\ldots\cdot\tau_m\}$$
the degree of $\pi$. The degree of $\pi$ depends only on its cycle type and is
sub-multiplicative: 
$$|\pi\rho|\leq |\pi|\cdot|\rho|\quad\quad\forall\pi\rho\in S_n.$$
Since $\pi\mapsto\sgn(\pi)=(-1)^{|\pi|}$
is a homomorphism, the {\em degree defect}
$$a(\pi_1,\pi_2,\ldots,\pi_t):=\half(|\pi_1|+\ldots+|\pi_t|-|\pi_1\cdots\pi_t|)$$
is a non-negative integer. 

For any subgroup $H\subset S_n$ and an $H$-stable subset $B\subset[n]:=\{1,\ldots,n\}$,
$H\backslash B$ denotes the orbit space for the induced action. For instance,
the degree can be expressed by $|\pi|=n-|\orb{\pi}{[n]}|$. 

For $\pi,\rho\in S_n$ the {\em graph defect} $g(\pi,\rho):
\orb{\pi,\rho}{[n]}\lra \IQ$ is defined by 
$$g(\pi,\rho)(B)=
\half\left(|B|+2-|\orb{\pi}{B}|-|\orb{\rho}{B}|-|\orb{\pi\rho}{B}|\right).$$

\begin{lemma}\label{lemma:orientedsurface}--- 
The graph defect $g$ takes value in the non-negative integers.
\end{lemma}

\begin{proof} Treating each orbit separately, we may assume that $B=[n]$, i.e.\
that the group $\<\pi,\rho\>$ acts transitively on $[n]$. We will identify the
defect $g$ as the genus of an oriented closed compact surface $C$ so that 
the claim becomes obvious.

Let $\sigma:=(\pi\rho)^{-1}$. Let $C$ be constructed as follows: Take 
$\{1,\ldots,n\}$ as the set of vertices. For all $i\in\{1,\ldots,n\}$ and 
$g\in\{\pi,\rho,\sigma\}$ add an oriented edge $e_{i,g}$ from $i$ to $g(i)$. 
For each $i$ glue in a (black) triangle along the edges $e_{i,\sigma}$, 
$e_{\sigma(i),\rho}$, and $e_{\rho\sigma(i),\pi}$. Finally, for each 
$g\in\{\pi,\rho,\sigma\}$ and each orbit $B^\prime$ of $g$ glue in a
(white) $|B^\prime|$-gon along the edges $e_{i,g}$, $e_{g(i),g}$,\ldots,
$e_{g^{|B^\prime|-1}(i),g}$ for some element $i\in B^\prime$. 
(If $|B^\prime|$ is $1$ or $2$, then the $|B^\prime|$-gon is rather degenerate: 
a disc with $1$ respectively $2$ edges.) Every edge bounds one white and one black
polygon. The resulting CW-complex is a connected compact oriented surface $C$ with
Euler characteristic
\begin{eqnarray*}
2-2g(C)=\chi(C)&=&\#\text{vertices}-\#\text{edges}+\#\text{polygons}
\\&=&n-3n+(a+b+c+n),
\end{eqnarray*}
where $a$, $b$, and $c$ are the numbers of black polygons corresponding to
the orbits of $\pi$, $\rho$, and $\sigma$, respectively. Hence
$g(C)=\half\left(n+2-a-b-c\right)=g(\pi,\rho)$.
\end{proof}

As in the case of the degree defect the graph defect allows a natural 
multivariant extension: for $\pi_1,\ldots\pi_t\in S_n$ let 
$$g(\pi_1,\ldots,\pi_t)(B)=\half
\left(|B|+2-\sum_j|\orb{\pi_j}{B}|-|\orb{\pi_1\cdots\pi_t}{B}|\right).$$
for $B\in \orb{\pi_1,\ldots,\pi_t}{[n]}$.


\subsection{} Now consider the graded vector space
\begin{equation}
{A\{S_n\}}:=\bigoplus_{\pi\in S_n}
A^{\tensor\orb{\pi}{[n]}}\cdot \pi.
\end{equation}
Here the grading of an element $a\cdot \pi$ is $|a\cdot \pi|:=|a|$.
The symmetric group $S_n$ acts on $A\{S_n\}$: the action of $\sigma\in S_n$ on
$[n]$ induces a bijection 
$$\sigma:\orb\pi{[n]}\to\orb{\sigma\pi\sigma^{-1}}{[n]},\quad x \mapsto \sigma x$$
for each $\pi$ and hence an isomorphism
\begin{equation}\label{eq:sigmaoperation}
\tilde\sigma:{A\{S_n\}}\lra{A\{S_n\}},\quad a\pi\mapsto 
\sigma^*(a)\sigma\pi\sigma^{-1}.
\end{equation}
Let
$$A^{[n]}:=(A\{S_n\})^{S_n}$$
be the subspace of invariants. 

\begin{example}--- For $n=3$ we get
$$A\{S_3\}=A^{\tensor 3}\id\oplus A^{\tensor 2}(12)\oplus
A^{\tensor 2}(13)\oplus A^{\tensor 2}(23)\oplus
A(123)\oplus A(132),$$
where we agree to give the orbits the lexicographical order with respect to 
$1<2<3$. Symmetrising, we find
$$A^{[3]}\isom S^3A\oplus (A\tensor A)\oplus A.$$
\end{example}


\subsection{}
Let $\kv(A):=S^*(A\tensor t^{-1}\IQ[t^{-1}])$ be the bosonic
Fock space modelled on the graded vector space $A$. Then $\kv(A)$ is bigraded
by degree and weight, where an element $a\tensor t^{-m}\in A\tensor t^{-m}$
is given degree $|a|$ and weight $m$. The component of $\kv(A)$ of constant weight $n$
is the graded vector space
$$\kv(A)_n\isom\bigoplus_{||\alpha||=n}\bigotimes_{i}S^{\alpha_i}A,$$
where $\alpha=(1^{\alpha_1},2^{\alpha_2}\cdots)$ runs through all partitions of
$n$ and $||\alpha||:=\sum_i i\alpha_i$.

If we think of $\IQ$ as a Frobenius algebra of degree $d=0$, then
$\IQ\{S_n\}=\IQ[S_n]$ is the group ring and $\kv(\IQ)=\IQ[p_1,p_2,\ldots]$ with
$p_m=t^{-m}$ is the ring of symmetric functions. The following map generalizes
the classical characteristic map $\IQ[S_n]\to\IQ[p_1,p_2,\ldots]$.

Let $f:\{1,\ldots,N\}\to \orb{\pi}{[n]}$ be an enumeration of the orbits of 
$\pi\in S_n$, and let
$\ell_i:=|f(i)|$ denote the length of the $i$-th orbit. Then define
$$\Phi':A^{\tensor N}\lra \kv(A), a_1\tensor\cdots\tensor a_N\mapsto \frac{1}{n!}
(a_1\tensor t^{-\ell_1})\cdots(a_N\tensor t^{-\ell_N}),$$
and let  
\begin{equation}\label{eq:PhiDefinition}
\Phi:\bigoplus_{n\geq 0}A\{S_n\}\lra \kv(A)
\end{equation}
be given on the summand $A^{\tensor\orb{\pi}{[n]}}\pi$ by the composition
$$A^{\tensor \orb{\pi}{[n]}}\xra{\ f^{-1}\ }A^{\tensor N}\xra{\Phi'} \kv(A).$$

\begin{proposition}\label{pr:PhiProposition}---
$\Phi$ induces an isomorphism of graded vector spaces 
$$A^{[n]}\lra \kv(A)_n.$$
\end{proposition}

\begin{proof} The map $\Phi$ is surjective and invariant under the 
$S_n$-action, so that its restriction to $A^{[n]}$ is also surjective.
Moreover, an invariant vector $v$ in $A^{[n]}$ is determined by its components
$v_\pi$ in $A^{\tensor\orb{\pi}{[n]}}$, where $\pi$ runs through a system of 
representatives for all possible cycle types, i.e.\ all partitions of $n$,
and each component $v_n$ is symmetric with respect to exchanging the values
corresponding to orbits of the same length, and conversely. This shows that 
both vector spaces have the same dimensions.
\end{proof}

Our next goal is to put a ring structure on $A\{S_n\}$ such that $A^{[n]}$ 
becomes a commutative subring.


\subsection{}
Any inclusion $H\subset K$ of subgroups of $S_n$ leads to a
surjection $H\backslash[n]\epimorph K\backslash[n]$ of orbit spaces
and hence to maps
\begin{equation}
f^{H,K}:A^{\tensor H\backslash[n]}\lra A^{\tensor K\backslash[n]}\quad
\text{and}\quad f_{K,H}:A^{\tensor K\backslash[n]}\lra A^{\tensor
H\backslash[n]}.
\end{equation}
If $H=\<\pi\>$ is the cyclic subgroup generated by a permutation $\pi$, we
omit the brackets $\<-\>$ in the notation. For $\pi,\rho\in S_n$ define
\begin{equation}
\begin{array}{c}
m_{\pi,\rho}:A^{\tensor\orb{\pi}{[n]}}\tensor A^{\tensor\orb{\rho}{[n]}}
\lra A^{\tensor\orb{\pi\rho}{[n]}}\\[1.5ex]
m_{\pi,\rho}(a\tensor b):=f_{\<\pi,\rho\>,\pi\rho}
\left(f^{\pi,\<\pi,\rho\>}(a)\cdot f^{\rho,\<\pi,\rho\>}(b)
\cdot e^{g(\pi,\rho)}\right)
\end{array}
\end{equation}
where $g(\pi,\rho):\orb{\pi,\rho}{[n]}\to\IN_0$ is the graph defect defined
in \ref{subsec:graphdefect}.

\begin{proposition}--- The product $A\{S_n\}\times A\{S_n\}\xra{\ \cdot\ }A\{S_n\}$
defined by
$$a\pi\ \cdot\  b\rho:=m_{\pi,\rho}(a\tensor b)\pi\rho$$
is associative, $S_n$-equivariant, and homogeneous of degree $nd$.
\end{proposition}

\begin{proof} For any $\pi$, $\rho$, $\sigma\in S_n$ consider the following
diagram of orbit spaces:
{\small
$$
\begin{array}{ccccccccc}
\orb\pi{[n]}&&\orb\rho{[n]}&&\orb{\pi\rho}{[n]}&&\orb\sigma{[n]}
&&\orb{\pi\rho\sigma}{[n]}\\[,5ex]
&\searrow&\downarrow&\swarrow&&\searrow&\downarrow&\swarrow\\[.5ex]
&&\orb{\pi,\rho}{[n]}&&\star&&\orb{\pi\rho,\sigma}{[n]}\\[.5ex]
&&&\searrow&&\swarrow\\
&&&&\orb{\pi,\rho,\sigma}{[n]}
\end{array}
$$
}The arrows of type $\searrow$ and $\downarrow$ will correspond to 
ring homomorphisms where\-as arrows of type $\swarrow$ will contravariantly
induce module homomorphisms. By definition, we have
$(a\pi\cdot b\rho)\cdot\,c\sigma =Y\pi\rho\sigma$
with 
$$Y=f_{\<\pi\rho,\sigma\>,\pi\rho\sigma}
\left(f^{\pi\rho,\<\pi\rho,\sigma\>}f_{\<\pi,\rho\>,\pi\rho}(Y')
\cdot f^{\sigma,\<\pi\rho,\sigma\>}(c)\cdot e^{g(\pi\rho,\sigma)}\right).
$$
and
$$Y'=f^{\pi,\<\pi,\rho\>}(a)\cdot f^{\rho,\<\pi,\rho,\>}(b)\cdot e^{g(\pi,\rho)}.$$
Lemma \ref{lemma:diamond} applies to the central diamond $\star$ in the diagram 
above: let $\nu_\star:\orb{\pi,\rho,\sigma}{[n]}\to \IN_0$ be defined by
$$\nu_\star(B)=1-|\orb{\pi,\rho}{B}|-|\orb{\pi\rho,\sigma}{B}|
+|\orb{\pi\rho}{B}|.$$
Then by \ref{lemma:diamond}:
\begin{eqnarray*}
\lefteqn{f^{\pi\rho,\<\pi\rho,\sigma\>}f_{\<\pi,\rho\>,\pi\rho}(Y')}
\hspace{1em}\\
&=&f_{\<\pi,\rho,\sigma\>,\<\pi\rho,\sigma\>}
(f^{\<\pi,\rho\>,\<\pi,\rho,\sigma\>}(Y')\cdot e^{\nu_\star})\\
&=&f_{\<\pi,\rho,\sigma\>,\<\pi\rho,\sigma\>}
\left(
f^{\pi,\<\pi,\rho,\sigma\>}(a)\cdot f^{\rho,\<\pi,\rho,\sigma\>}(b)\cdot
f^{\<\pi,\rho\>,\<\pi,\rho,\sigma\>}(e^{g(\pi,\rho)})\cdot e^{\nu_\star}
\right).
\end{eqnarray*}
We apply the projection formula and get
\begin{equation}
Y=f_{\<\pi,\rho,\sigma\>,\pi\rho\sigma}\left(
f^{\pi,\<\pi,\rho,\sigma\>}(a)\cdot 
f^{\rho,\<\pi,\rho,\sigma\>}(b)\cdot
f^{\sigma,\<\pi,\rho,\sigma\>}(c)\cdot
Y''\right),
\end{equation}
with
$$Y''=f^{\<\pi,\rho\>,\<\pi,\rho,\sigma\>}(e^{g(\pi,\rho)})
\cdot e^{\nu_\star}\cdot
f^{\<\pi\rho,\sigma\>,\<\pi,\rho,\sigma\>}(e^{g(\pi\rho,\sigma)})
=:e^h.
$$
In this expression $h:\orb{\pi,\rho,\sigma}{[n]}\lra\IN_0$ is the function
\begin{eqnarray*}
h(B)&=&\sum_{B^\prime\in\orb{\pi,\rho}{B}}g(\pi,\rho)(B^\prime)+
\nu_\star(B)+\sum_{B^\dprime\in\orb{\pi\rho,\sigma}{B}}g(\pi\rho,\sigma)(B^\dprime)
\\
&=&\half\left(|B|+2|\orb{\pi,\rho}{B}|-|\orb{\pi}{B}|-|\orb{\rho}{B}|
-|\orb{\pi\rho}{B}|\right)\\
&&+\left(1-|\orb{\pi,\rho}{B}|-|\orb{\pi\rho,\sigma}{B}|-|\orb{\pi\rho}{B}|
\right)\\
&&+\half\left(|B|+2|\orb{\pi\rho,\sigma}{B}|-|\orb{\pi\rho}{B}|
-|\orb{\sigma}{B}|-|\orb{\pi\rho\sigma}{B}|\right)\\
&=&|B|+1-\half\left(|\orb{\pi}{B}|+|\orb{\rho}{B}|+|\orb{\sigma}{B}|+
|\orb{\pi\rho\sigma}{B}|\right)\\
&=&g(\pi,\rho,\sigma)(B).
\end{eqnarray*}
Summing up, we have
$$
Y=f_{\<\pi,\rho,\sigma\>,\pi\rho\sigma}\left(
f^{\pi,\<\pi,\rho,\sigma\>}(a)\cdot 
f^{\rho,\<\pi,\rho,\sigma\>}(b)\cdot
f^{\sigma,\<\pi,\rho,\sigma\>}(c)\cdot
e^{g(\pi,\rho,\sigma)}\right).
$$
The same symmetric expression arises if we compute $a\chi_\pi\cdot (b\chi_\rho\cdot
c\chi_\sigma)$ in the same way. This proves associativity.

The following terms contribute to the degree of $|a\pi\cdot b\rho|$:
the degrees of $a$, $b$ and $e^{g(\pi,\rho)}$ (given by Lemma \ref{edegree}), 
the degrees of $f^{\pi,\<\pi,\rho\>}$, 
$f^{\rho,\<\pi,\rho\>}$ and $f_{\<\pi,\rho\>,\pi\rho}$, and the degrees of the
two multiplications in the definition of $m_{\pi,\rho}$.
The total balance is:
\begin{eqnarray*}\lefteqn{|a\pi\cdot b\rho|}\quad\\
&=&|a|+|b|+|f^{\pi,\<\pi,\rho\>}|+|f^{\rho,\<\pi,\rho\>}|
+|f_{\<\pi,\rho\>,\pi\rho}|+|e^{g(\pi,\rho)}|+2d|\orb{\pi,\rho}{[n]}|\\[1ex]
&=&|a|+|b|+d(|\orb{\pi\rho}{[n]}|-|\orb{\pi,\rho}{[n]}|)
+d(|\orb{\pi}{[n]}|-|\orb{\pi,\rho}{[n]}|)\\[1ex]
&&+d(|\orb{\rho}{[n]}|-|\orb{\pi,\rho}{[n]}|)
+ 2d\sum_{B\in 
\orb{\pi,\rho}{[n]}}g(\pi,\rho)(B)-d|\orb{\pi,\rho}{[n]}|\\[-1.5ex]
&&+2|\orb{\pi,\rho}{[n]}|\\
&=&|a\pi|+|b\rho|+nd
\end{eqnarray*}
\end{proof}

\begin{proposition}\label{pr:centerrelation}--- For any two homogeneous elements
$a\pi, b\rho\in A\{S_n\}$ the following (non)commutativity relation holds:
$$a\pi\cdot b\rho=(-1)^{|a|\cdot|b|}\pi^*(b)\pi\rho\pi^{-1}\cdot a\pi=
(-1)^{|a|\cdot|b|}\tilde\pi(b\rho)\cdot a\pi.$$
\end{proposition}

\begin{proof} Let $\rho^\prime:=\pi\rho\pi^{-1}$. The following diagram of
orbit spaces
$$
\begin{array}{ccccccc}
&&\orb{\rho}{[n]}&
\xra{\quad\pi\quad}&\orb{\rho^\prime}{[n]}\\[1ex]
&&\big\downarrow&&\big\downarrow\\[1ex]
\orb{\pi}{[n]}&\longrightarrow&\orb{\pi,\rho}{[n]}&\xra{\quad =\quad}
&\orb{\rho^\prime,\pi}{[n]}&\longleftarrow&\orb{\pi\rho}{[n]}
\end{array}
$$
commutes. Hence 
$f^{\rho^\prime,\<\pi,\rho^{\prime}\>}(\pi^*(b))=f^{\rho,\<\pi,\rho\>}(b)$.
Moreover, $g(\pi,\rho)=g(\rho^{\prime},\pi)$. It follows that
\begin{eqnarray*}m_{\pi,\rho}(a\tensor b)
&=&f_{\<\pi,\rho\>,\pi\rho}(f^{\pi,\<\pi,\rho\>}(a)\cdot f^{\rho,\<\pi,\rho\>}(b)
\cdot e^{g(\pi,\rho)})\\
&=&f_{\<\rho^{\prime},\pi\>,\pi\rho}(f^{\pi,\<\pi,\rho\>}(a)\cdot 
f^{\rho^{\prime},\<\pi,\rho^{\prime}\>}(\pi^*(b))\cdot 
e^{g(\rho^{\prime},\pi)})\\
&=&(-1)^{|a|\cdot|b|}
f_{\<\rho^{\prime},\pi\>,\pi\rho} (f^{\rho^{\prime},\<\pi,\rho^{\prime}\>}
(\pi^*(b))\cdot f^{\pi,\<\pi,\rho\>}(a) \cdot e^{g(\rho^{\prime},\pi)})\\
&=&(-1)^{|a|\cdot|b|}m_{\rho^{\prime},\pi}(\pi^*(b)\tensor a)
\end{eqnarray*}
\end{proof}

\begin{proposition}--- $A^{[n]}$ is a subring of the centre of $A\{S_n\}$.
\end{proposition}

\begin{proof} This is a consequence of Proposition \ref{pr:centerrelation}.
\end{proof} 

\begin{proposition}--- Let $\Frobint:{A\{S_n\}}\to \IQ$ be defined by 
$$\Frobint(a\pi):=\left\{\begin{array}{ll}\Frobint(a)&\text{if }\pi=\id,\\0&
\text{else},\end{array}\right.$$
where $\Frobint$ on the right hand side is the integral on $A^{\tensor [n]}$. 
The restriction of this integral to $A^{[n]}$ defines the structure of a graded
Frobenius algebra of degree $nd$ on $A^{[n]}$.
\end{proposition}

\begin{proof} The integral induces an $S_n$-invariant bilinear form on $A\{S_n\}$. 
The only non-trivial pairings are of type
$$A^{\tensor \orb{\pi}{[n]}}   \tensor 
A^{\tensor\orb{\pi^{-1}}{[n]}}\lra k.$$
As $\pi$ and $\pi^{-1}$ have the same orbit spaces, this map is the composition
$$A^{\tensor\orb{\pi}{[n]}}\tensor A^{\tensor\orb{\pi}{[n]}}\xra{\ \cdot\ }
A^{\tensor\orb{\pi}{[n]}}\xra{\varphi_*}A^{\tensor n}\xra{\Frobint} k.$$
But $\Frobint\circ\varphi_*=\Frobint$, the integral on $A^{\tensor\orb{\pi}{[n]}}$.
This shows that $\Frobint$ induces a non-degenerate pairing on $A\{S_n\}$.
Since the pairing is invariant, its restriction to $A^{[n]}$ is also 
non-degenerate.
\end{proof}

\begin{example}--- To illustrate the multiplication in $A\{S_n\}$ we take up
the preceding example $n=3$. One finds:
$$(\alpha\tensor \beta)(12)\cdot 
(\gamma\tensor\delta)(13)=\alpha\beta\gamma\delta(132),$$
$$(\alpha\tensor\beta)(12)\cdot
(\gamma\tensor \delta)(12)=(-1)^{|\beta|\cdot 
|\gamma|}\Delta_*(\alpha\gamma)\tensor(\beta\delta) \id,$$
$$\alpha(123)\cdot \beta(123)=(\alpha\beta e)(132),$$
$$\alpha(123)\cdot \beta (132)=\Delta_*(\alpha\beta)\id.$$
\end{example}


\section{Hilbert schemes}


\subsection{}
If one applies the construction of the 
previous section to the cohomology ring of a smooth projective surface
with trivial canonical divisor one obtains the cohomology ring of the
Hilbert schemes. In order to do so we need to shift the
degrees of the cohomology rings by the dimension of the corresponding
manifold. 

Let $X$ be a smooth projective surface over the complex numbers. Let
$X^{[n]}:=\Hilb^n(X)$ denote $n$-th Hilbert scheme, i.e.\ the moduli
space that represents the functor
$$S\mapsto\{Z\subset S\times X\text{ closed subscheme }|p:Z\to S
\text{ flat, finite of degree }n\}.$$
Then $X$ is again projective \cite{Grothendieck} and smooth \cite{Fogarty}
of dimension $2n$. The following is the main theorem of this paper:

\begin{theorem}\label{th:MainTheorem}--- 
Let $X$ be a smooth projective surface with numerically
trivial canonical divisor. Then there is a canonical isomorphism of graded
rings
$$(H^*(X;\IQ)[2])^{[n]}\xra{\quad\isom\quad}H^*(X^{[n]};\IQ)[2n].$$
\end{theorem}


\subsection{}
The proof will be given in the next section. We recall some results on
Hilbert schemes and their cohomology which will be used in the proof.
We also take the opportunity to change some sign conventions from 
\cite{Lehn} to adapt our formulae to sign rules which are standard
in the literature about vertex algebras.

Let $H:=H^*(X;\IQ)[2]$ and $\IH_n:=H^*(X^{[n]},\IQ)[2n]$. Then 
$\IH:=\bigoplus_{n\geq 0}\IH_n$ is bigraded by the (shifted) cohomological
degree and the conformal weight $n$. The shift has the effect of centring the
middle degrees of the cohomology groups at zero. There is a distinguished
element $\vacuum\in H^0(X^{[0]};\IQ)\subset \IH$, the vacuum.

As compact complex manifolds, $X$ and the Hilbert schemes 
$X^{[n]}$ have fundamental classes in the top degree homology groups. We
denote the evaluation of cohomology classes on the fundamental class by
$\int_{[X]}$ and $\int_{[X^{[n]}]}$, respectively. Now give $H$ and 
$\IH_{n}$ the structure of graded Frobenius algebras by setting
$\Frobint(a):=-\int_{[X]}a$ for $a\in H$ and 
$\Frobint(a):=(-1)^n\int_{[X^{[n]}]}a$ for $a\in \IH_n$.
Note that the two conventions agree on 
$X=X^{[1]}$.

\subsection{}
For $n\in \IN$, Nakajima defines incidence varieties
$$Z_n:=\{(\xi,x,\xi')|\, \xi\subset \xi',
|\xi'|-|\xi|=nx\}$$
in $X^{[\ell]}\times X\times X^{[\ell+n]}$, and operators
$\Nakop_{-n}: H\to\End_\IQ(\IH)$, 
$$\Nakop_{-n}(\alpha)(y):= PD^{-1}(pr_{3*}((pr_2^*(\alpha)\cup pr_1^*(y))\cap 
[Z_n])),$$
where $y\in \IH_\ell$, and $PD$ denotes Poincar\'e duality.
Moreover, we define 
$$\Nakop_{n}(\alpha):=\Nakop_{-n}(\alpha)^\dagger,$$
where for an endomorphism $f\in \End_\IQ(\IH)$, $f^\dagger$ denotes the adjoint
endomorphism with respect to the pairing on $\IH$ given by $\Frobint$. Finally,
to simplify  statements, we let $\Nakop_0(\alpha)=0$. The operators 
$\Nakop_{-n}(\alpha)$ are bihomogeneous of bidegree $(n,|\alpha|)$.

\begin{theorem}--- For $n,m\in \IZ$ and $a,b\in H$, the following oscillator
relation holds: $[\Nakop_n(a),\Nakop_m(b)]=n\cdot\delta_{n,-m}\cdot 
\Frobint(ab)\cdot\id_\IH$.
\end{theorem}

Here and in the following all elements in $H$ and $\IH$ and, accordingly,
endomorphisms of $\IH$ are graded by their cohomological degree. The commutator
in the theorem and those in the remainder of this paper have to be taken in the 
graded sense. The theorem is due to Nakajima \cite{Nakajima} and Grojnowski
\cite{Grojnowski}, the coefficient $n$ 
on the right hand side was determined by Ellingsrud and Str\o{}mme \cite{E-S3}.
Combining the theorem with the calculation
of the Betti numbers of $X^{[n]}$ due to G\"ottsche one obtains by formal 
arguments:

\begin{theorem}\label{th:Nak}--- There is an isomorphism of graded vector spaces 
$$\Psi:\kv(H)\lra \IH,\quad (a_1t^{-n_1})\cdots 
(a_st^{-n_s})\mapsto\Nakop_{-n_1}(a_1)\cdots\Nakop_{-n_s}(a_s)\vacuum.$$
\end{theorem}

Recall that $\kv(H)=S^*(H\tensor t^{-1}\IQ[t^{-1}])$.
Next, let $\Xi_n\subset X^{[n]}\times X$ denote the universal subscheme, and
let $p:\Xi_n\to X^{[n]}$ and $q:\Xi_n\to X$ denote the two projections. Then 
for any element $\alpha\in H$, let 
$$\alpha^{[n]}:=p_*(ch(\ko_{\Xi_n})\cdot q^*(td(X)\cdot \alpha)).$$
We denote by $\alpha^{[\bullet]}\in\End_\IQ(\IH)$ the linear operator which on 
$\IH_n$ is multiplication by the class $\alpha^{[n]}$. Of particular importance
is the class $1^{[n]}=ch(\ko^{[n]})$, the Chern character of the tautological
sheaf $\ko^{[n]}:=p_*(\ko_{\Xi_n})$. The homogeneous component of degree 2 of the
operator $1^{[\bullet]}$ will be denoted by $\partial$ (see \cite{Lehn}). The following 
proposition was first proved by the first named author \cite[Thm 4.2]{Lehn} for the
special case $\alpha=ch(F)$, $F$ a locally free sheaf on $X$, and then extended
to the general case by Li, Qin and Wang \cite{LiQinWang1}:

\begin{theorem}\label{th:L-LQW}--- For all $\alpha,y \in H$ one has
$$[\alpha^{[\bullet]},\Nakop_{-1}(y)]=\exp(\rm{ad}(\partial))\Nakop_{-1}
(\alpha\cdot y).$$
\end{theorem}

The theorem requires the computation of the (higher order) commutator of 
$\partial$ and $\Nakop_{-n}$. Let $\Delta_*:H\to H\tensor H$ denote the 
adjoint of the multiplication map. Then for any two integers $n,m\in \IZ$ and
an element $a\in H$ we obtain an operator $\Nakop_n\Nakop_m(\Delta_*(a))$.
Define
$$L_n(a):=\frac{1}{2}\sum_{\nu\in\IZ}\,:\Nakop_{\nu}\Nakop_{n-\nu}:\Delta_*(a),
$$
where $:-:$ denotes the normal ordered product of two operators (i.e.\ 
\mbox{$:\Nakop_{\nu}\Nakop_{n-\nu}:$}$\ =\Nakop_{\nu}\Nakop_{n-\nu}$ if
$n-\nu>0$ and
$\Nakop_{n-\nu}\Nakop_{\nu}$ otherwise). The definition of the operators 
$L_n(a)$ is a twisted version of the definition of the standard conformal
structure for a free super boson \cite{Kac}. These operators satisfy the 
relations
$$[L_n(a),L_m(b)]=(n-m)L_{n+m}(ab)+\delta_{n,-m}\frac{n^3-n}{12}T(e_Hab),$$
where $e_H$ is the Euler class of $H$, which is $-c_2(X)$ by our conventions,
and $$[L_n(a),\Nakop_m(b)]=-m\Nakop_{n+m}(ab).$$
The following theorem is the main result of \cite{Lehn}:

\begin{theorem}\label{th:LMain}--- For $n\in \IN$ and $a\in H$ one has
$$[\partial, \Nakop_{-n}(a)]=(-n) L_{-n}(a)+\binom{-n}{2}\Nakop_{-n}(Ka),$$
where $K\in H$ is the class of the canonical divisor of $X$.
\end{theorem}

Iterated application of the theorem leads to the identity
$$\textstyle{\frac{(-1)^{n}}{n!}}
{\rm ad}^n([\partial,\Nakop_{-1}(1)])(\Nakop_{-1}(a))=\Nakop_{-n-1}(a),$$
from which it is clear that
\begin{equation}\label{eq:inductionstart}
\IH=\partial\IH+\Nakop_{-1}(H)\IH.
\end{equation}
This opens the path for many induction arguments on $\IH$, since $\partial$
increases the cohomological degree, and $\Nakop_{-1}(H)$ the weight.

Theorems \ref{th:Nak}, \ref{th:L-LQW} and \ref{th:LMain} together give a 
complete explicit description of the ring structure of $\IH_n$. It is proved
in \cite{LiQinWang1} that the elements $a^{[n]}$, $a\in H$, form a 
set of generators. Although the relations among these generators are explicitly
given in the sense that the theorems above provide implementable algorithms to
compute all products etc., the problem remains to give a description of the
resulting ring in terms of workable generators and relations. In this sense
this paper is related to \cite{LiQinWang1} as our paper \cite{LehnSorger} to
the last section of \cite{Lehn}. The following proposition follows formally
from Theorem \ref{th:LMain}:

\begin{proposition}\label{pr:firstchernclass}--- 
Assume that $K=0$.  Let $n=n_1+\ldots+n_s$
and let $a_1,\ldots,a_s\in H$ be homogeneous classes. Then
\begin{eqnarray*}
\lefteqn{c_1(\ko^{[n]})\cdot \Nakop_{-n_1}(a_1)\cdots\Nakop_{-n_s}(a_s)\vacuum}
\hspace{0em}\\
&=&-\sum_{i<j}\eps_{ij}\,n_in_j\Nakop_{-n_i-n_j}(a_ia_j)
\Nakop_{-n_1}(a_1)\cdots\widehat\Nakop_{-n_i}\cdots\widehat\Nakop_{-n_j}
\cdots\Nakop_{-n_s}(a_s)\vacuum\\
&&-\half\sum_{i}\eps_i\sum_{n'+n''=n_i}n_i\Nakop_{-n'}\Nakop_{-n''}\Delta_*(a_i)
\Nakop_{-n_1}(a_1)\cdots\widehat\Nakop_{-n_i}\cdots\Nakop_{-n_s}(a_s)\vacuum,
\end{eqnarray*}
where the $\eps$'s account for the signs which result from commuting the 
$a_k$'s and $\widehat\Nakop$ indicates operators that are omitted.
\end{proposition}

\begin{proof} By definition, multiplication by $c_1(\ko^{[n]})$ is the same
as applying the operator $\partial$. Now move the operator $\partial$ as far
to the right as possible thus introducing commutators. The claim of the 
proposition then follows from the following facts:
\begin{enumerate}
\item $\partial \vacuum=0$,
\item $[[\partial,\Nakop_{-n}(a)],\Nakop_{-m}(b)]=-nm\Nakop_{-n-m}(ab)$, 
\item 
$[\partial,\Nakop_{-n}(a)]\vacuum=-\half\sum_{n'+n''=n}
n\Nakop_{-n'}\Nakop_{-n''}(\Delta_*(a))$.
\end{enumerate}
Here (1) holds degree reasons, and (2) and (3) follow directly from 
\ref{th:LMain}.\end{proof}


\section{Proof of Theorem \ref{th:MainTheorem}}\label{sc:TheProof}


\subsection{}
We keep the notations of the previous sections. In the following $X$ will
always be a smooth projective surface with numerically
trivial canonical divisor. $H=H^*(X;\IQ)[2]$ is a graded
Frobenius algebra of degree $d=2$ in the sense of the first section.
There are isomorphisms of bigraded vector spaces
$$\Gamma:\bigoplus_{n\geq 0}H^{[n]}\xra{\quad\Phi\quad}\kv(H)\xra{\quad\Psi\quad}
\bigoplus_{n\geq 0} H^*(X^{[n]};\IQ)[2n]=:\IH.$$
Here $\Phi$ is the map of Proposition \ref{pr:PhiProposition} and $\Psi$ is 
Nakajima's isomorphism in Theorem \ref{th:Nak}. The proof of Theorem 
\ref{th:MainTheorem} will be obtained by carefully identifying elements in 
and operators on these three spaces and translating bits of information about
one ring into information about the other rings. By the very definition of
$\Psi$ the operator $a\tensor t^{-m}$ on $\kv(H)$ corresponds to Nakajima's
operator $\Nakop_{-m}(a)$ on $\IH$. Analogously, there are operators 
$\gothr_{-m}(a)$ on $\bigoplus_{n\geq 0}H^{[n]}$. We will only need the 
operator $\gothr_{-1}(a)$ and give an ad hoc definition:

For any element $y\in H\{S_n\}$ and $a\in H$ let $y\tensor a\in H\{S_{n+1}\}$
denote the element that is obtained by adding the trivial cycle $(n+1)$ to each
permutation $\pi\in S_n$ and tensoring the coefficient of $\pi$ with $a$.
Moreover, for all $n\in \IN_0$ let $P:H\{S_n\}\to H^{[n]}$ be the 
symmetrisation operator 
$$P(y)=\frac{1}{n!}\sum_{\sigma\in S_n}\widetilde\sigma(y).$$
Then
$$\gothr_{-1}(a)(y)=(n+1)\cdot(-1)^{|a|\cdot|y|}P(y\tensor a).$$
Note that this defines in fact a linear map $\gothr_{-1}(a):H\{S_n\}\to 
H^{[n+1]}$. But for symmetric $y\in H^{[n]}$,
we may simplify the symmetrisation operator:
$$(n+1)\cdot P(y\tensor a)=\sum_{i=1}^{n+1}(i\,n+1)^{\sim}(y\tensor a),$$
where $(i\,n+1)^\sim$ is the action of the transposition $(i\,n+1)$ as given by
(\ref{eq:sigmaoperation}), and where by abuse of notation we allow $(i\,n+1)$ to 
denote the identity permutation if $i=n+1$.

Let $\Gamma:=\Psi\Phi$.
With the given definition of $\gothr_{-1}$ one has
$$\Gamma(\gothr_{-1}(a)(y))=\Psi((a\tensor t^{-1})(\Phi(y)))
=\Nakop_{-1}(a)(\Gamma(y)).$$

\subsection{} Let $\eps_n:=\sum_{\pi}\sgn(\pi)\pi\in H^{[n]}$ denote the 
alternating character. Recall that $\ko^{[n]}:=p_*\ko_{\Xi_n}$ denotes the
tautological sheaf on $X^{[n]}$.

\begin{proposition}--- $\Gamma(\eps_n)=c(\ko^{[n]})$.
\end{proposition}

\begin{proof} Recall that the number of permutations $\pi$ of
a given cycle type 
$\lambda=(\lambda_1,\lambda_2,\ldots)=(1^{\alpha_1}2^{\alpha_2}\cdots)$
is given by
$$C_\lambda=\frac{n!}{\prod_i \alpha_i! i^{\alpha_i}},$$
and that its signature is given by 
$\sgn(\pi)=(-1)^{\sum_{j=1}^s(\lambda_j-1)}=:\sgn(\lambda).$
Therefore
\begin{eqnarray*}\Gamma(\eps_n)
&=&\frac{1}{n!}\sum_{\lambda\in\text{partitions of }n}
C_\lambda\sgn(\lambda)
\Nakop_{-\lambda_1}(1)\cdots\Nakop_{-\lambda_{s}}(1)\vacuum\\
&=&\sum_{||\alpha||=n}\prod_{i}\frac{1}{\alpha_i!}
\left(\frac{(-1)^{i-1}\Nakop_{-i}(1)}{i}\right)^{\alpha_i}\vacuum\\
&=&\exp\left(\sum_{i}\frac{(-1)^{i-1}}{i}\Nakop_{-i}(1)\right)_n\vacuum,
\end{eqnarray*}
where the index $n$ in the last line means that we take the component of
weight $n$ only. By \cite[Theorem 4.6]{Lehn}, the last expression equals
the total Chern class of the tautological sheaf $\ko^{[n]}$.
\end{proof}

If we pick only the component of degree $2$, we get $c_1(\ko^{[n]})$ on the
right hand side and
$$\eps_{n,2}:=-\sum_{\text{all transpositions}}\tau$$
on the left hand side.

\begin{proposition}\label{pr:alternatingone}--- 
$\Gamma(\eps_{n,2}\cdot y)=c_1(\ko^{[n]})\cdot 
\Gamma (y)$ for all $y\in H^{[n]}$.
\end{proposition}

\begin{proof} This proposition is an adaptation of a result of Goulden 
\cite[Proposition 3.1.]{Goulden} to our situation. Without loss of generality
we may assume that $y$ is of the form $y=P(a\pi)$. Moreover, $\eps_{n,2}\cdot y=
P(\eps_{n,2}\cdot a\pi)$. Therefore,
we may fix a permutation with a disjoint cycle decomposition $\pi=z_1\cdots z_s$ 
with a certain ordering and assume that $y=(a_1\tensor\cdots\tensor a_s) \pi$. 
(Of course, then $y$ is no longer contained in $H^{[n]}$, but recall that $\Phi$
is defined on the larger ring $H\{S_n\}$.) Then by definition,
$$\Gamma(y)=\Nakop_{-\ell_1}(a_1)\cdots\Nakop_{-\ell_s}(a_s)\vacuum.$$

Now let $\tau$ be a single transposition, say $\tau=(i\,j)$. We distinguish two cases
according to whether $i$ and $j$ are contained in the same $\pi$-orbit or not. 
We analyse the affect of multiplying $\tau$ and $\pi$, and take the sum over all
transpositions afterwards.

{\em 1.\ case}: $i$ and $j$ are contained in different cycles of $\pi$, say
$$\pi=(i\,x_2\cdots x_\ell)(j\,z_2\cdots z_m)\cdots.$$
Then $\tau\pi=(i\,x_2\cdots x_\ell\,j\,z_2\cdots z_m)\cdots$, i.e.\ the 
multiplication with $\tau$ merges the two orbits. The genus defect 
$g(\tau,\pi)$ vanishes. Hence the multiplication map
\begin{equation}\label{eq:gouldeneq1}
m_{\tau,\pi}(1\tensor -):H^{\tensor\orb{\pi}{[n]}}\lra 
H^{\tensor\orb{\tau\pi}{[n]}}
\end{equation}
is essentially given by multiplying the coefficients corresponding to the
two orbits $B':=\{i,x_2,\cdots,x_\ell\}$ and $B'':=\{j,z_2,\cdots, z_m\}$.
Assuming that $i$ and $j$ are contained in the $\tilde\imath$'s and $\tilde\jmath$'s
orbit of $\pi$, this consideration shows
\begin{equation}\label{eq:Gouldenfirstpart}
\Gamma(\tau\cdot y)=\eps^\prime\Nakop_{-\ell_{\tilde\imath}-\ell_{\tilde\jmath}}
(a_{\tilde\imath}a_{\tilde\jmath})\Nakop_{-\ell_1}(a_1)\cdots
\widehat\Nakop_{-\ell_{\tilde\imath}}\cdots\widehat\Nakop_{-\ell_{\tilde\jmath}}\cdots
\Nakop_{-\ell_s}(a_s),
\end{equation}
where $\eps^\prime$ is the sign arising from the permutation of the $a_k$'s.
If $\tau$ runs through all transpositions, there are $|B'|\cdot |B''|=
\ell_{\tilde\imath}\ell_{\tilde\jmath}$ possibilities to hit the orbits $B'$ and $B''$. 
Thus the right hand side in (\ref{eq:Gouldenfirstpart}) occurs with 
multiplicity $\ell_{\tilde\imath}\ell_{\tilde\jmath}$. Up to the sign, this yields
the first term on the right hand side of Proposition \ref{pr:firstchernclass}.

{2.\ case}: $i$ and $j$ are contained in the same cycle of $\pi$, say
$$\pi=(i\,x_2\cdots x_{m'}\,j\,z_2\cdots z_{m''})\cdots.$$
Then $\tau\pi=(i\,x_2\cdots x_{m'})(j\,z_2\cdots z_{m''})\cdots$, i.e.\ the given
cycle is split into two smaller cycles. Again the genus defect vanishes. The
multiplication (\ref{eq:gouldeneq1}) is essentially given by the 
comultiplication $H\to H\tensor H$, where the two factors on 
the right hand side correspond to the two new orbits. Hence if the cycle
$(i\,x_2\cdots x_{m'}\,j\,z_2\cdots z_{m''})$ is, say, $z_h$, then 
\begin{equation}\label{eq:Gouldensecondpart}
\Gamma(\tau\cdot y)=\eps^{\prime\prime}\Nakop_{-m'}\Nakop_{-m''}(\Delta_*(a_h))
\Nakop_{-\ell_1}(a_1)\cdots\widehat\Nakop_{-\ell_h}\cdots\Nakop_{-\ell_s}(a_s),
\end{equation}
where again $\eps^{\prime\prime}$ is the sign arising from the permutation of
the $a_k$'s. There are precisely $\ell_h$ choices of (ordered!) pairs $(i,j)$ from
the cycle $z_h$ such that the cycle splits into two cycles of lengths $m'$ and 
$m''$. Again up to the sign, this corresponds to the second term on the right hand
side of Proposition \ref{pr:firstchernclass}. Note that the factor $\half$ 
arises since transpositions are unordered.

Summing up we see, that multiplication by $\eps_{n,2}$ has the same effect
-- via $\Gamma$ -- as multiplication by 
$c_1(\ko^{[n]})$ as described by Proposition \ref{pr:firstchernclass}.
\end{proof}


\subsection{}
It follows from \cite[Theorem 4.2.]{Lehn} that
\begin{equation}\label{eq:chernclassrelation}
c(\ko^{[n+1]})\Nakop_{-1}(a)-\Nakop_{-1}(a)c(\ko^{[n]})=[\partial, 
\Nakop_{-1}(a)]c(\ko^{[n]})
\end{equation}
as operators on $\IH_{n}$. We must prove that the corresponding assertion holds
for the alternating character. 

\begin{proposition}\label{pr:alternatingrelation}
--- The following identity of operators on $H^{[n]}$ holds:
$$\eps_{n+1}\gothr_{-1}(a)-\gothr_{-1}(a)\eps_{n}=
(\eps_{n+1,2}\gothr_{-1}(a)-\gothr_{-1}(a)\eps_{n,2})\eps_{n}$$
\end{proposition}

\begin{proof} The embedding $S_n\to H\{S_n\}$ preserves products provided
that $|\pi\rho|=|\pi|\cdot|\rho|$. Therefore, there are identities
$$\eps_{n+1,2}-\iota(\eps_{n,2})=\sum_{i=1}^n (i\ n+1)\quad\text{and}\quad 
\eps_{n+1}-\iota(\eps_n)=\sum_{i=1}^n (i\ n+1)\cdot \iota(\eps_n).$$
The proposition follows from this by a simple calculation (see the proof of 
\cite[Proposition 5.1.]{LehnSorger}.
\end{proof}

\begin{proposition}\label{pr:alternatingall}--- The following identity holds
for all $y\in H^{[n]}$:
$$\Gamma(\eps_{n}\cdot y)=c(\ko^{[n]})\cdot\Gamma(y).$$
\end{proposition}

\begin{proof} The assertion follows from (\ref{eq:inductionstart}),
(\ref{eq:chernclassrelation}) and Proposition \ref{pr:alternatingrelation}
by induction on cohomological degree and weight. 
The calculation itself is identical to the proof 
of \cite[Proposition 5.3.]{LehnSorger}.
\end{proof}


\subsection{}
Let $TH^{[n]}\subset H^{[n]}$ and $T\IH_n\subset\IH_n$ denote the tautological
rings, i.e.\ the subalgebras generated by the components $\eps_{n,2k}$ of the
alternating character and  the components $c(\ko^{[n]})_{k}$ of the tautological
bundle, respectively. What we have proved so far can be rephrased as follows:
$\Gamma$ maps $TH^{[n]}$ to $T\IH_n$ and, more precisely, $TH^{[n]}\to T\IH_n$
is an isomorphism of rings, and $H^{[n]}\to \IH_n$ is an isomorphism of modules
over the tautological rings. It still remains to show that $\Gamma$ is an
isomorphism of rings. To see this we show that the ring structure of either ring
is completely determined by the module structure over the tautological subring.
The key point here is, that certain operators satisfy what Li, Qin and Wang
\cite{LiQinWang1} call the transfer property: for the classes $a^{[n]}$ it
follows directly from Theorem \ref{th:L-LQW} that
\begin{equation}
[a^{[\bullet]},\Nakop_{-1}(b)]=[1^{[\bullet]},\Nakop_{-1}(ab)]
\end{equation}
We proceed in three steps. In section \ref{subsec:schritt1} we identify the
elements $\Gamma^{-1}(a^{[n]})\in H^{[n]}$; in section \ref{subsec:schritt2}
we show that the elements $\Gamma^{-1}(a^{[n]})$ have the transfer property;
and, finally, in section \ref{subsec:schritt3} we use the transfer property
to complete the proof by induction.


\subsection{}\label{subsec:schritt1}
For $a\in H$ consider the sum $\sum_{n\geq 0}a^{[n]}$ (in the formal
completion of $\IH$ with respect to the filtration by conformal weight). 
We put $\gothp:=\gothp_{-1}(1)$ for convenience and get
\begin{eqnarray*}
\sum_{n\geq 0}a^{[n]}&=&a^{[\bullet]}\exp(\gothp)\vacuum\\[-2ex]
&=&\sum_{n\geq 0}\frac{1}{n!}\sum_{k=0}^{n-1}\gothp^{n-k-1}[a^{[\bullet]},\gothp]
\gothp^k\vacuum,\\[.5ex]
&&\text{since $a^{[\bullet]}\vacuum=0$,}\\[.5ex]
&=&\sum_{n\geq 0}\sum_{k=0}^{n-1}\frac{1}{n!}\gothp^{n-k-1}\sum_{s=0}^k\binom{k}{s}
\gothp^{k-s}(-\ad\gothp)^s([a^{\bullet]},\gothp])\vacuum\\
&=&\sum_{m=0}^\infty\frac{\gothp^m}{m!}\cdot\sum_{s=1}^{\infty}\frac{1}{s!}
(-\ad\gothp)^{s-1}([a^{[\bullet]},\gothp])\vacuum.
\end{eqnarray*}
Inserting 
$$[a^{[\bullet]},\gothp]=\exp(\ad\partial)\gothp_{-1}(a)$$
we get
\begin{equation}\label{eq:schritt1:1}
\sum_{n\geq 0}a^{[n]}=\exp(\gothp)\sum_{s=1}^{\infty}
\sum_{k=0}^\infty\frac{(-\ad\gothp)^{s-1}}{s!}
\frac{(\ad\partial)^k}{k!}(\gothp_1(a))\vacuum.
\end{equation}
Let $\alpha$ be a partition of length $|\alpha|$ and let $\Delta_*:H\to 
H^{\tensor|\alpha|}$ be the diagonal map as before. For $a\in H$ we define
the operator
$$\gothp_{-\alpha}(a):=\prod_{i\geq 1}(\gothp_{-i})^{\alpha_i}(\Delta_*(a)).$$

\begin{proposition}\label{pr:thecoefficients}--- 
For each partition $\alpha$ there are rational numbers
$c_\alpha',c_\alpha''$ such that for $c_\alpha:=c_\alpha'+c_\alpha''e\in H$
the following identity holds:
\begin{equation}\label{eq:schritt1:2}
\sum_{s=1}^{\infty}\sum_{k=0}^\infty\frac{(-\ad\gothp)^{s-1}}{s!}
\frac{(\ad\partial)^k}{k!}(\gothp_1(a))\vacuum=
\sum_\alpha\frac{1}{\|\alpha\|!}\gothp_{-\alpha}(ac_\alpha)\vacuum.
\end{equation}
\end{proposition}

Pondering the left hand side the accuracy of this proposition becomes ma\-ni\-fest,
even more so as we are only claiming something about the structure of the right
hand side and nothing about the coefficients themselves. However, we find it painful to
formally prove the assertion and know no simpler way than doing some vertex algebra
calculus. We postpone the proof to subsection \ref{subsec:stupidproof}.

By (\ref{eq:schritt1:1}) and (\ref{eq:schritt1:2}) we have
\begin{equation}\label{eq:schritt1:3}
a^{[n]}=\sum_{||\alpha||\leq n}\frac{\gothp^{n-\|\alpha\|}}{(n-\|\alpha\|)!}
\frac{\gothp_{-\alpha}(c_\alpha a)}{\|\alpha\|!}\vacuum.
\end{equation}


\subsection{}\label{subsec:schritt2}
For any partition $\alpha=(1^{\alpha_1}2^{\alpha_2}\cdots)$ of 
$m=\|\alpha\|=\sum_ii\alpha_i$ of length $|\alpha|=\sum_i\alpha_i$ choose a 
permutation $\pi\in S_m$ of cycle type $\alpha$. As before, let
$\Delta_*:H\to H^{\tensor|\alpha|}$ be the map adjoint to the multiplication. 
For all $n\in \IN_0$ let
$$B_{\alpha}(u)_n:=\binom{n}{m}P(\Delta_*(u)\pi\tensor 1^{\tensor n-m})
\in H^{[n]},$$
where $P$ is the symmetrisation operator. One checks that
\begin{equation}\label{eq:schritt2:1}
\Gamma(B_{\alpha}(u)_n)=\frac{\gothp^{n-\|\alpha\|}}{(n-\|\alpha\|)!}
\cdot \frac{\gothp_{-\alpha}(u)}{\|\alpha\|!}\vacuum.
\end{equation}
Furthermore, let $B_\alpha(u)_\bullet$ denote the operator on
$\bigoplus_{n\geq 0}H^{[n]}$ which is multiplication by the elements
$B_\alpha(u)_n$.

\begin{proposition}\label{pr:transfer}--- The elements $B_\alpha(u)_n$ satisfy 
$$
[B_\alpha(u)_{\bullet},\gothr_{-1}(b)]=[B_\alpha(1)_{\bullet},\gothr_{-1}(ub)].
$$
In particular, $\Gamma^{-1}(a^{[n]})$ has the transfer property
$$[\Gamma^{-1}(a^{[\bullet]}),\gothr_{-1}(b)]=[\Gamma^{-1}(1^{[\bullet]}),
\gothr_{-1}(ab)].$$
\end{proposition}

\begin{proof} Let $y\in H^{[n]}$ be given. We compute the terms 
$B_\alpha(u)_{n+1}\cdot\gothr_{-1}(b)(y)$ and 
$\gothr_{-1}(ub)(B_\alpha(1)_{n}\cdot y)$:
\begin{eqnarray*}
\lefteqn{(-1)^{|y|\cdot|b|}B_\alpha(u)_{n+1}\cdot\gothr_{-1}(b)(y)}
\hspace{4em}\\[1ex]
&=&
\textstyle{\binom{n+1}{m}}\,P((\Delta_*(u)\pi\tensor 1^{\tensor n+1-m})\cdot 
(n+1)P(y\tensor b))\\[2ex]
&=&\textstyle{\binom{n+1}{m}}\,P\left((\Delta_*(u)\pi\tensor 1^{\tensor n+1-m})\cdot 
\sum_{i=1}^{n+1}(i\,n+1)^\sim (y\tensor b)\right)\\[2ex]
&=&\textstyle{\binom{n+1}{m}}\,\sum_{i=1}^{n+1}
P\left(((i\,n+1)^{\sim}(\Delta_*(u)\pi\tensor 1^{\tensor n+1-m}))\cdot
(y\tensor b)\right)
\end{eqnarray*}
If in this sum the index $i$ takes values $m+1,\ldots,n+1$, then the
transposition $(i\,n+1)$ permutes the final 1's in $\Delta_*(u)\tensor 
1^{\tensor n+1-m}$ and thus has no effect. This part of the sum therefore
equals
\begin{eqnarray*}
\phantom{4em}&=&\textstyle{\binom{n+1}{m}(n+1-m)}
P((\Delta_*(u)\pi\tensor 1^{\tensor n-m}\cdot y)\tensor b)\\[2ex]
&=&\textstyle{\binom{n+1}{m}\frac{n+1-m}{n+1}}\ \gothr_{-1}(b)
(\Delta_*(u)\tensor 1^{\tensor n-m}\cdot y)\\[2ex]
&=&(-1)^{(|y|+|u|)|b|}\gothr_{-1}(b)(B_\alpha(u)_n\cdot y).
\end{eqnarray*}
This is the second half of the commutator. Adding up, we find
\begin{eqnarray*}
\lefteqn{[B_\alpha(u)_{\bullet},\gothr_{-1}(b)]}\hspace{2em}\\
&=&(-1)^{|b|\cdot|u|}\textstyle{\binom{n+1}{m}}\,\sum_{i=1}^mP\left((i\,n+1)^\sim
(\Delta_*(u)\tensor 1^{n+1-m}\cdot (y\tensor b)\right).
\end{eqnarray*}
There is no need to evaluate the right hand side. It suffices to note that 
conjugation by $(i\,n+1)$ has the effect of making the point $(n+1)\in[n+1]$
part of the orbits which support the diagonally embedded class $\Delta_*(u)$. By
the definition of the multiplication in $A\{S_{n+1}\}$, this leads to a 
contraction of $u$ and $b$. But the result of the multiplication does not
change if we replace the pair $(u,b)$ by $(1,ub)$ or $(ub,1)$. This is all we need.
The second claim follows from the first when combined with 
(\ref{eq:schritt1:3}) and (\ref{eq:schritt2:1}).
\end{proof}


\subsection{}\label{subsec:schritt3}
%
Now we are ready to finish the proof of Theorem \ref{th:MainTheorem}. 
We must prove that for all $a\in H$ and homogeneous $y\in H^{[n]}$ one has
\begin{equation}\label{eq:tobeshown}
\Gamma(\Gamma^{-1}(a^{[n]})\cdot y)=a^{[n]}\cdot \Gamma(y).
\end{equation}
Since $1^{[n]}=ch(\ko^{[n]})$ is contained in the tautological ring, we already 
know that
\begin{equation}\label{eq:alreadyknown}
\Gamma(\Gamma^{-1}(1^{[n]})\cdot y)=1^{[n]}\cdot \Gamma(y).
\end{equation}

We argue by induction on conformal weight and cohomological degree and 
assume that (\ref{eq:tobeshown}) is true for all conformal weights $<n$ 
and all degrees $<|y|$. We know that 
$$\IH_n=\partial\IH_{n}+\Nakop_{-1}(1)\IH_{n-1}\quad\text{and}\quad
H^{[n]}=\eps_{n,2}\cdot H^{[n]}+\gothr_{-1}(1)H^{[n-1]}.$$
It therefore suffices to consider the following two cases:

\noindent
{\em Case 1}: $y=\eps_{n,s}\cdot z$. Then
\begin{eqnarray*}
\Gamma(\Gamma^{-1}(a^{[n]})\cdot y)
&=&\Gamma(\Gamma^{-1}(a^{[n]})\cdot\eps_{n,2}\cdot z)
=\Gamma(\eps_{n,s}\cdot \Gamma^{-1}(a^{[n]})\cdot z)\\
&=&\Gamma(\eps_{n,2})\cdot \Gamma(\Gamma^{-1}(a^{[n]})\cdot z)
=\partial\cdot a^{[n]}\Gamma(z)\\
&&\text{ by induction,}\\
&=&a^{[n]}\partial\Gamma(z)=a^{[n]}\cdot\Gamma(\eps_{n,2}\cdot z)=
a^{[n]}\Gamma(y).
\end{eqnarray*}

\noindent
{\em Case 2}: $y=\gothr_{-1}(1)z$. Here we make use of the transfer property
(Proposition \ref{pr:transfer}) which explicitly says:
\begin{eqnarray}\nonumber
\lefteqn{\Gamma^{-1}(a^{[n]})\gothr_{-1}(1)-\gothr_{-1}(1)\Gamma^{-1}(a^{[n-1]})}
\hspace{4em}\\
\label{eq:doppelzeile}
&=&\Gamma^{-1}(1^{[n]})\gothr_{-1}(a)-\gothr_{-1}(a)\Gamma^{-1}(1^{[n-1]}).
\end{eqnarray}
Then
\begin{eqnarray*}
\Gamma(\Gamma^{-1}(a^{[n]})\cdot y)
&=&\Gamma(\Gamma^{-1}(a^{[n]})\cdot\gothr_{-1}(1)z)\\
&=&\Gamma(\gothr_{-1}(1)\Gamma^{-1}(a^{[n-1]})\cdot z)\\
&&+\Gamma(\Gamma^{-1}(1^{[n]})\gothr_{-1}(a)z)
-\Gamma(\gothr_{-1}(a)\Gamma^{-1}(1^{[n-1]})\cdot z)\\
&&\text{because of (\ref{eq:doppelzeile}),}\\
&=&\Nakop_{-1}(1)\Gamma(\Gamma^{-1}(a^{[n-1]})z)\\
&&+1^{[n]}\cdot\Gamma(\gothr_{-1}(a)\cdot z)
-\Nakop_{-1}(a)1^{[n-1]}\cdot\Gamma(z)\\
&&\text{because (\ref{eq:tobeshown}) holds true for $1^{[\nu]}\in T\IH$,}\\
&=&\Nakop_{-1}(1)a^{[n-1]}\Gamma(z)\\
&&+1^{[n]}\gothp_{-1}(a)1^{[n-1]}\Gamma(z)
-\Nakop_{-1}(a)1^{[n-1]}\cdot\Gamma(z)\\
&&\text{by induction,}\\
&=&a^{[n]}\Nakop_{-1}(1)\Gamma(z)\\
&&\text{because of the transfer property,}\\
&=&a^{[n]}\Gamma(\gothr_{-1}(1)z)=a^{[n]}\Gamma(y).
\end{eqnarray*}
This finishes the proof of the main theorem (up to the proof of Proposition 
\ref{pr:thecoefficients} in the next section).

\subsection{Proof of Proposition \ref{pr:thecoefficients}}
\label{subsec:stupidproof}

Let $\glf(\IH)\subset \End(\IH)[|z,z^{-1}|]$ denote the general linear field
algebra of $\IH$. Our basic fields are
$$\varphi(a)(z):=\sum_{n\in\IZ}\gothp_{n}(a)z^{-n-1}$$
for $a\in H$ and their derivatives
$$\varphi(a)^{(k)}(z):=\left(\frac{\partial}{\partial z}\right)^{k}\varphi(a)(z).$$
More generally, for any partition $\beta=1^{\beta_1}2^{\beta_2}\cdots$ let
$$\varphi_{\beta}(a)(z):=\quad
:\prod_i\left(\frac{\varphi^{(i-1)}}{i!}\right)^{\beta_i}\Delta_*(a):,$$
where as before $\Delta_*:H\to H^{\tensor\|\beta\|}$ is the map adjoint to
multiplication.

Then $\varphi_{\beta}(a)(z)$ is a field of conformal weight $\|\beta\|$. We 
recover the operators $\gothp(a)$ and $\partial$ as Fourier modes of
 $\varphi(a)(z)$ and $\frac{1}{3!}\varphi_{1^3}(1)(z)$:
$$\varphi(a)(z)=\ldots+\gothp(a)\cdot z^0+\ldots\quad\text{and}\quad\frac{1}{3!}
\varphi_{1^3}(1)(z)=\ldots+(-\partial)\cdot z^{-3}+\ldots$$
(For the latter fact see \cite{Frenkel-Wang}). The Wick Theorem 
(cf.\ \cite{Kac}[Thm.\ 3.3.]) applies to these fields
and yields the following OPE:
\begin{eqnarray*}
\varphi(1)(z)\cdot \varphi_\beta(b)(w)&\sim&
\sum_{j\geq 1}\beta_j\varphi_{\beta-j^1}(b)(w)\frac{1}{j!}
\frac{\partial^{j-1}}{\partial w^{j-1}}\frac{1}{(z-w)^2}\\
&\sim&\sum_{j\geq 1}\frac{\beta_j\varphi_{\beta-j^1}(b)(w)}{(z-w)^{j+1}},
\end{eqnarray*}
which implies the commutator relation (cf.\ \cite{Kac}[Thm.\ 2.3.])
$$-\ad\gothp\big(\varphi_\beta(b)(w)\big)=
\sum_{j\geq 1}(-w)^{-j-1}\beta_j\varphi_{\beta-j^1}(b)(w).$$
Here and in the following $\beta-j^1$ denotes the partition that equals $\beta$ 
with the number of $j$'s decreased by one, similarly for $\beta+j^1-k^2$ etc.

In a similar way, we compute the OPE for the fields 
$\frac{1}{3!}\varphi_{1^3}(1)(z)$ and $\varphi_\beta(b)(z)$. In this case, the Wick
Theorem gives several terms on the right hand side depending on whether we
contract $1$, $2$ or $3$ factors. Observe, however, that contracting $N$
factors introduces the power $e^{N-1}$ of the Euler class. As $e^2$ is zero,
only the following terms are left:
\begin{eqnarray*}
\lefteqn{\frac{1}{3!}\varphi_{1^3}(1)(z)\cdot\varphi_\beta(b)(w)\sim\frac{1}{2}
\sum_{j\geq 1}\beta_j:\varphi(z)^2\varphi_{\beta-j^1}(w):(b)\frac{1}{(z-w)^{j+1}}}
\hspace{4em}\\
&&+\sum_{1\leq k<\ell}\beta_k\beta_\ell:\varphi(z)\varphi_{\beta-k^1-\ell^1}(w):(be)
\frac{1}{(z-w)^{k+\ell+2}}\\
&&+\sum_{1\leq k}\binom{\beta_k}{2}:\varphi(z)\varphi_{\beta-k^2}(w):(be)
\frac{1}{(z-w)^{2k+2}}\\
&\lefteqn{\text{which by Taylor expansion yields:}}\\
&
\sim&\frac{1}{2}\sum_{j\geq 1}\sum_{k,\ell\geq 1}\beta_jk\ell\,
:\varphi_{\beta-j^1+k^1+\ell^1}(b)(w):\,\frac{1}{(z-w)^{j-k-\ell+3}}\\
&
&+\sum_{1\leq j}\sum_{1\leq k<\ell}\beta_k\beta_\ell j\,
:\varphi_{\beta+j^1-k^1-\ell^1}(be)(w):\,\frac{1}{(z-w)^{k+\ell-j+3}}\\
&
&+\sum_{1\leq j}\sum_{1\leq k}\binom{\beta_k}{2}j\,
:\varphi_{\beta+j^1-k^2}(be)(w):\,\frac{1}{(z-w)^{2k-j+3}}.
\end{eqnarray*}
As before we pass to the commutator relation for the Fourier mode $\partial$
of $\frac{1}{3!}\varphi_{1^3}(1)(z)$ and get:
\begin{eqnarray*}
(-\ad\partial)\big(\varphi_\beta(b)(w)\big)&=&\sum_{m=0}^2\binom{2}{m}w^m\Big\{
\frac{1}{2}\sum_{\substack{j=k+\ell-m\\1\leq k,\ell}}
\beta_jk\ell\varphi_{\beta-j^1+k^1+\ell^1}(b)(w)\\
&&+\sum_{\substack{j=k+\ell+m\\1\leq k<\ell}}\beta_k\beta_\ell j
\varphi_{\beta+j^1-k^1-\ell^1}(be)(w)\\
&&+\sum_{\substack{j=2k+m\\1\leq k}}\binom{\beta_k}{2}j\varphi_{\beta+j^1-k^2}(be)(w)
\Big\}
\end{eqnarray*}
These commutator relations can be expressed more elegantly as follows:
The map
\begin{eqnarray*}
\mu:H[z,z^{-1}][t_1,t_2,\ldots]&\xra{\quad\isom\quad}&\IQ[z,z^{-1}]\langle
\varphi_{\beta}(a)(z)\rangle_{a,\beta}\subset\glf(\IH)\\
at_1^{\beta_1}t_2^{\beta_2}\cdots t_s^{\beta_s}
&\mapsto&\varphi_\beta(a)(z)z^{\|\beta\|}
\end{eqnarray*}
is a $\IQ[z,z^{-1}]$-linear isomorphism onto the submodule generated by the 
fields $\varphi_{\beta}(a)(z)$. With respect to this identification, the operators
$\ad\gothp$ and $\ad\partial$ as calculated above can be written as 
$$-\ad\gothp=z^{-1}D_1\quad\text{and}\quad\ad\partial=-D_2$$
with differential operators
$$
D_1:=\sum_{j\geq 1}(-1)^{j-1}\frac{\partial}{\partial t_j}$$
and
$$
D_2:=\frac{1}{2}\sum_{k+\ell\geq j}
\binom{2}{k+\ell-j}k\ell t_kt_\ell\frac{\partial}{\partial t_j}
+\frac{e}{2}\sum_{j\geq k+\ell}\binom{2}{j-k-\ell}jt_j
\frac{\partial}{\partial t_k}\frac{\partial}{\partial t_\ell}.
$$
To compute the left hand side of (\ref{eq:schritt1:2}) we calculate 
in $H[z,z^{-1}][t_1,t_2,\ldots]$. Recall that $at_1\mapsto \varphi(a)(z)\cdot z$ and
that $\gothp_{-1}(a)$ is the coefficient of $z$ in this field. Then 
$$f:=\exp(-D_2)(t_1)$$
is a power series in $\IQ[e][|t_1,t_2,\ldots|]$
and we need the coefficient of $z$ in the field corresponding
to $a\cdot\sum_{s\geq 1}\frac{1}{s!}z^{1-s}D_1^{s-1}f$ when applied to the vacuum, or 
equivalently, the term
$$\sum_{s\geq 1}\mbox{Coeff}\left(z^{s},\mu\Big( \frac{a}{s!}D_1^{s-1}f\Big)
\vacuum\right).$$
To apply the field $\varphi$ or normal ordered products of its 
derivatives to the vacuum simply requires to
throw away all Fourier modes $\gothp_{m}$, $m\geq 0$, i.e. we evaluate
$$\tilde\mu:t_i\mapsto \sum_{n>0}\frac{\gothp_{-n}}{n}\binom{n}{i}z^{n-i}\in
\End(\IQ[\gothp_{-1},\gothp_{-2},\ldots])[z].$$
Now we are done: For 
$$
\sum_{s\geq 1}\frac{1}{s!}\mbox{Coeff}(z^{s},\tilde\mu(D_1^{s-1}f))
\in \IQ[e][|\gothp_{-1},\gothp_{-2},\ldots|]
$$
is a power series in $\gothp_{-\nu}$'s and hence can be expressed as
$$\sum_{s\geq 1}\frac{1}{s!}\mbox{Coeff}(z^{s},\tilde\mu(D_1^{s-1}f))
=
\sum_{\alpha=(1^{\alpha_1}\cdots r^{\alpha_r})}\frac{c_\alpha}{\|\alpha\|!}
\gothp_{-1}^{\alpha_1}\gothp_{-2}^{\alpha_2}\cdots\gothp_{-r}^{\alpha_r}$$
for suitable coefficients $c_\alpha\in\IQ+\IQ e\subset H$, so that in total
$$
\sum_{s=1}^{\infty}\sum_{k=0}^\infty\frac{(-\ad\gothp)^{s-1}}{s!}
\frac{(\ad\partial)^k}{k!}(\gothp_{-1}(a))\vacuum=
\sum_{\alpha}\frac{1}{\|\alpha\|!}\gothp_{-\alpha}(c_\alpha\cdot a)\cdot\vacuum.$$

\begin{remark}--- Numerical evidence suggests the simple expression
$$c_\alpha=\frac{(-1)^{\|\alpha\|-|\alpha|}}{\prod_i\alpha_i!}\left(1+
\frac{|||\alpha|||-1}{24}e\right),$$
where $|||\alpha|||:=\sum_i\binom{i+1}{2}\alpha_i$. 
\end{remark}

\bibliographystyle{plain}

\parindent0mm

\end{document}